\documentclass[11pt,a4paper]{amsart}

\usepackage{amsmath,mathrsfs,amssymb}

\newtheorem{theorem}{Theorem}[section]
\newtheorem{lemma}[theorem]{Lemma}
\newtheorem{proposition}[theorem]{Proposition}
\newtheorem{corollary}[theorem]{Corollary}

\newtheorem{question}[theorem]{Question}

\theoremstyle{definition}

\theoremstyle{plain}

\newcommand{\FamilyF}{\mathscr{F}}

\newcommand{\ol}[1]{\overline{#1}}
\newcommand{\loopauto}[2]{\ensuremath{\hat\Gamma_{#1}({#2})}}

\begin{document}
    \title[The Loop Problem]{The Loop Problem \\ for Monoids and Semigroups}

\maketitle

\begin{center}
    Mark Kambites \\

    \medskip

    School of Mathematics, \  University of Manchester \\
    Manchester M60 1QD, \  England. \\

    \medskip

    \texttt{Mark.Kambites@manchester.ac.uk} \\

\end{center}

\begin{abstract}
We propose a way of associating to each finitely generated monoid or
semigroup a formal language, called its \textit{loop problem}. In the
case of a group, the loop problem is essentially the same as the
\textit{word problem} in the sense of combinatorial group theory.
Like the word problem for groups, the loop problem is regular if and only
if the monoid is finite. We study also the case in which the loop problem is
\textit{context-free}, showing that a celebrated group-theoretic result of
Muller and Schupp extends to describe completely simple semigroups with
context-free loop problems. We consider also right cancellative monoids,
establishing connections between the loop problem and the
structural theory of these semigroups by showing that the syntactic monoid of the
loop problem is the \textit{inverse hull} of the monoid.
\end{abstract}

\bigskip

\section{Introduction}\label{sec_intro}

One of the most productive and successful approaches to finitely
generated groups is the study of the language of all words over a
given finite generating set which represent
the identity element. Many
of the structural properties of a group are reflected in the
language-theoretic properties of this language. For example,
the language is regular, context-free,
recursive or recursively enumerable exactly if the group is finite,
virtually free, embeddable in a simple subgroup of a finitely presented
group, or embeddable in a finitely presented group respectively
\cite{Anisimov71,Muller83,Boone74,Higman61}.

In a group, words $u$ and $v$ over the generators represent the same
element exactly if $u v^{-1}$ represents the identity. It follows that
deciding membership for the language of words representing the identity
is algorithmically equivalent to the \textit{word problem} in the sense of
universal algebra \cite[Section~30]{Graetzer68}, that is, the problem of deciding whether two given words in the
generators represent the same element of the group. In particular, the
language encodes the
multiplication, and hence the entire structure, of the group. For this
reason, the language of words representing the identity is usually termed the
\textit{word problem} of the group.

By contrast, in a more general finitely generated monoid, there is no
equivalence between the language of words representing the identity, and
the universal algebraic word problem. One need only consider the class of
monoids obtained by adjoining identities to semigroups, to see that radically different monoids can give
rise to the same language.
It is perhaps for this reason that recent research in combinatorial semigroup
theory has, with a few exceptions \cite{Duncan04,KambitesHypCS},
eschewed direct interaction with formal language theory. Although automata
appear with increasing frequency (as in the emerging theory of
\textit{automatic semigroups} \cite{Hudson96, Campbell01}), they are
typically of a ``multi-tape'' form which does not admit a direct
language-theoretic interpretation. However, in view of the numerous and deep 
interactions between combinatorial group theory and formal languages, and
also between finite semigroup theory and formal languages (see, for example,
\cite{Pin86}), it would be surprising if the theory of formal languages had
no role to play in the combinatorial theory of finitely generated semigroups
and monoids.

The purpose of this article is to propose and begin the study of a new way
of associating to each finitely generated monoid or semigroup a formal
language, which we call the \textit{loop problem}. Like the word problem
for groups, this language encapsulates the entire structure of the semigroup,
at least in the simplistic sense that it \textit{uniquely determines} the semigroup. A less trivial
question is that of whether there is a correspondence between ``natural''
structural properties of a semigroup or monoid, and the ``natural''
language-theoretic properties of its loop problem, and it is this question in
which we are primarily interested.

In addition to this introduction, this paper comprises six
sections. In Section~\ref{sec_autos}, we briefly revise some elementary
automata theory. Section~\ref{sec_loop} introduces the loop problem of
a finitely generated monoid or semigroup, while Section~\ref{sec_basic}
establishes some of its basic properties, including the extent to which
it is invariant under change of generators.

Section~\ref{sec_groupcs} considers the  case in which
the monoid is a group, describing the precise relationship between the
loop problem and the word problem. More generally, we study the important class
of finitely generated completely simple semigroups, establishing
a natural relationship between the loop problems of such semigroups and
the loop problems (and word problems) of their maximal subgroups.

Section~\ref{sec_correspond} begins the
study of the interaction between natural structural properties of semigroups,
and language-theoretic properties of their loop problems. We begin with a
semigroup-theoretic analogue of a basic result of Anisimov
\cite{Anisimov71} from combinatorial group theory --- a monoid is finite
if and only if its loop problem is
regular. We consider also the classes of context-free languages; we apply a
result from Section~\ref{sec_groupcs} to show that a celebrated theorem of
Muller and Schupp \cite{Muller83} from combinatorial group theory extends
to characterise completely simple semigroups with context-free loop problems,
and ask what can be said more generally about the class of semigroups with
this property.

Finally, in Section~\ref{sec_canc}, we consider the loop problem in the
important special case of \textit{right cancellative} monoids. It transpires that there is a simple and elegant relationship between the
between the loop problem and the established structural theory of
these monoids \cite{Clifford61,Rees48}. In particular, the syntactic monoid of the
loop problem is the \textit{inverse hull} of the monoid.

\section{Automata}\label{sec_autos}

In this paper, a central role will be played by infinite automata over
finite alphabets. It is well-known, if not so well documented, that
much of the classical theory of finite automata extends easily to the
infinite case, provided sufficient care is taken with the definitions. In
this section, we briefly introduce some fundamental
definitions and results of automata theory, partly so as to
make the paper accessible to the reader not familiar with classical automata
theory, and partly to establish the necessary foundations in the infinite
case. Proofs are not given; they can be found in the finite case in the
standard textbooks on formal language theory \cite{Hopcroft69,Pin86}, while
the appropriate adjustments for the infinite case are left as an exercise
for the reader.

Let $X$ be a finite set of symbols or \textit{letters}, called an
\textit{alphabet}. A \textit{word} over $X$ is a finite sequence of
letters from $X$; the \textit{empty word}, with no symbols, is denoted $\epsilon$.
The collections of all words over $X$ forms a monoid (with identity $\epsilon$)
under the operation of concatenation; it is called the \textit{free monoid}
on $X$ and denoted $X^*$. The non-empty words comprise a subsemigroup of $X^*$,
called the \textit{free semigroup} on $X$ and denoted $X^+$. It is readily
verified that free monoids and semigroups satisfy the usual universal
property defining free objects in the categories of monoids and semigroups
respectively. A \textit{language} over $X$ is a collection of words over
$X$, that is, a subset of the free monoid $X^*$.

An \textit{automaton} over $X$ consists of a
directed graph with edges labelled by elements of $X^*$, with a
distinguished \textit{start vertex} and a set of distinguished
\textit{terminal vertices}. The vertices and edges of an automaton are
sometimes called \textit{states} and \textit{transitions} respectively.
The labelling of edges extends naturally,
via the multiplication in $X^*$, to a labelling of paths by words in
$X^*$. The \textit{language accepted} by the automaton is the set of
all words in $X^*$ which label paths from the initial vertex to a
terminal vertex.

The automaton is called \textit{finite} if the vertex and edge sets are
finite; a language accepted by a finite automaton is called
a \textit{regular language}. Many textbooks on automata focus
exclusively on the finite/regular case, and use the terms ``automaton''
and ``finite automaton'' interchangeable. In this article, we are
equally interested in the infinite case, and an ``automaton'' may be
infinite unless explicitly stated otherwise.

An automaton is called \textit{deterministic} if its edges are labelled
by elements of $X$, and for each state $q$ and
each letter $x \in X$, there is at most one edge leaving $q$ with label
$x$. An automaton is called \textit{trim} if for every vertex $q$, there
is a path from the initial state to $q$, and a path from $q$ to some
terminal state.
Two automata are \textit{isomorphic} if there is an isomorphism between
their underlying graphs which preserves edge labels, the start state and
the set of terminal states.

Let $Q$ be the state set of the automaton. For each word $w \in X^*$, one
can define a binary relation
$$\sigma_w = \lbrace (p, q) \mid \text{there is a path from $p$ to $q$ labelled $w$} \rbrace \subseteq Q \times Q.$$
The set of all relations of the form $\sigma_w$ constitutes a subsemigroup
of the monoid of all binary relations on $Q$; the multiplication is given
by $\sigma_u \sigma_v = \sigma_{uv}$.
The relation $\sigma_\epsilon$ 
is easily seen to be act as an identity for this semigroup (even though it
may not be the identity function on $Q$), so in fact the semigroup
is a monoid, called the \textit{transition monoid} of the automaton.
In the case that the automaton is deterministic, the relations $\sigma_w$ are
all partial functions, $\sigma_\epsilon$ is actually the identity
function on $Q$, and the transition monoid is generated by the maps $\sigma_x$
for letters $x \in X$.

The \textit{cone} of a vertex $q$ is the language of all words which label
paths from $q$ to some terminal vertex. Two vertices are called
\textit{equivalent} if they have the same cone.

\begin{proposition}\label{prop_minauto}
Let $L$ be a language. Then there exists a unique (up to isomorphism)
deterministic trim automaton $M(L)$ accepting $L$, with the property that
no two vertices are equivalent.
\end{proposition}

The automaton given by Proposition~\ref{prop_minauto} is called the
\textit{minimal automaton} of the language $L$. It
is a final object in the category of all deterministic trim automata
accepting $L$ with an appropriate notion of morphism. An easy corollary is
that $L$ is regular exactly if $M(L)$ is finite. Moreover, if $L$ is regular
then $M(L)$ has strictly
fewer states than any other deterministic automaton accepting $L$; many
authors concerned only with finite automata take this property as the definition
of the minimal automaton, although even in the finite case the morphism
viewpoint is arguably more helpful for understanding the theory. In the
infinite case, of course, the cardinality of the state set is a wholly
insufficient to characterize the minimal automaton.

Let $L$ be a language over an alphabet $X$. We say that two words
$u, v \in X^*$ are \textit{syntactically equivalent} (with respect
to $L$), and write
$u \equiv_L v$, if for every pair
of words $x, y \in L$, we have $xuy \in L$ if and only if $xvy \in L$.
Thus, two words are syntactically equivalent if one can substitute one
for the other in any word, without affecting membership of the language
$L$. We write $[u]_L$ for the syntactic equivalence class of a word $u$.
It is readily verified that the relation of syntactic equivalence is
a congruence on the free monoid $X^*$, so that the syntactic equivalence
classes form a monoid with multiplication well-defined by $[u]_L [v]_L = [uv]_L$.
This is called the
\textit{syntactic monoid} of $L$ and denoted $S(L)$. The following proposition
relates syntactic monoids to minimal automata.
\begin{proposition}\label{prop_syntrans}
Let $L$ be a language. Then the syntactic monoid $S(L)$ is isomorphic to
the transition monoid of the minimal automaton $M(L)$ via the map
$[u]_L \mapsto \sigma_u$.
\end{proposition}

We shall also need the notion of a \textit{rational transduction}. We provide
here a brief definition; for a detailed introduction, see \cite{Berstel79}.
Let $X$ and $Y$ be finite alphabets. A \textit{finite transducer} from $X$
to $Y$ is a finite directed graph with edges labelled by elements of the direct
product $X^* \times Y^*$, and with a distinguished initial vertex and a
set of distinguished terminal vertices. Just as with automata, the
labelling of edges extends, via the multiplication in the direct product
monoid $X^* \times Y^*$, to a labelling of paths. The \textit{relation
accepted} by the automaton is the set of all pairs in $X^* \times Y^*$ which
label paths from the initial vertex to a terminal vertex. A relation between
free monoids $X$ and $Y$ recognised by a finite transducer is called a
\textit{rational transduction} from $X$ to $Y$.

If $\rho \subseteq X^* \times Y^*$ is a relation and $L \subseteq X^*$ is
a language, then the \textit{image} of $L$ under $\rho$ is the language
$$L \rho \ = \ \lbrace v \mid (u, v) \in \rho \text{ for some } u \in L \rbrace \ \subseteq \ Y^*.$$
If $\rho$ is a rational transduction then we say that $L \rho$ is
\textit{a rational transduction of} $L$ (via $\rho$). Two languages are
called \textit{rationally equivalent} if each is a rational transduction
of the other; rational equivalence is, as the name implies, an equivalence
relation on the class of all languages. The \textit{inverse} of $\rho$ is the relation
$$\rho^{-1} \ = \ \lbrace (u, v) \mid (v, u) \in \rho \rbrace \subseteq Y^* \times X^*.$$
It is easily seen that if $\rho$ is a rational transduction then so is
$\rho^{-1}$.

We shall also need another operation on languages. If
$L \subseteq X^*$ is a language then the \textit{Kleene closure} of $L$
is the submonoid of $X^*$ generated by $L$, that is, the set of all words
of the form $w_1 \dots w_n$ for some integer $n \geq 0$ and words
$w_1, \dots w_n \in L$. The Kleene closure of $L$ is usually denoted $L^*$; note that if $L = X$ is in fact the
alphabet viewed as a subset of the free monoid it generates, then the Kleene
closure of $L$ is
the entire free monoid on $X$, so the notation $X^*$ is unambiguous.

\section{The Loop Problem}\label{sec_loop}

In this section, we introduce the languages which are the main object of
study in this paper. We assume familiarity with the basic terminology
and notation of semigroup theory
\cite{Howie95} and combinatorial group theory \cite{Lyndon77,Magnus76}. We define
the loop problem first for monoids, and then more generally for
semigroups.

Let $M$ be a monoid. By a \textit{choice of (monoid) generators} for
$M$ we mean a surjective morphism $\sigma : X^* \to M$ from a free monoid
onto $M$. The choice of generators is called \textit{finite} if $X$ is
finite. A monoid which admits a finite choice of generators is called
\textit{finitely generated}. Since our aim is to connect semigroup theory
with the theory of formal languages (over finite alphabets) we are
primarily interested in finitely generated monoids, although many of our
results will also hold for infinite choices of generators.

Recall that the \textit{(right) Cayley graph} $\Gamma_\sigma(M)$ of $M$ with
respect to $\sigma$ is a directed graph, possibly with multiple edges and
loops, with edges labelled by elements of $X$.
Its vertices are the elements of $M$, and it has an edge from $a \in M$ to
$b \in M$ labelled $x \in X$ exactly if $a (x \sigma) = b$ in the monoid $M$.

Now let $\ol{X} = \lbrace \ol{x} \mid x \in X \rbrace$
be a set of formal inverses for the generators in $X$, and let
$\hat{X} = X \cup \ol{X}$. We extend the map $x \to \ol{x}$ to an
involution on $\hat{X}^*$ by defining $\ol{\ol{x}} = x$ for all $x \in X$,
and $\ol{x_1 \dots x_n} = \ol{x_n} \dots \ol{x_1}$ for all
$x_1, \dots, x_n \in \hat{X}$. (In fact, $\hat{X}^*$ together with the
unary operation
$x \to \ol{x}$ is the \textit{free monoid with involution} on $X$ \cite[Section~6.1]{Lawson98}.)

The \textit{(right) loop automaton} $\loopauto{\sigma}{M}$ of $M$ with
respect to $X$ is obtained from the Cayley graph $\Gamma_\sigma(M)$ by
adding for each edge labelled $x$ an inverse edge, in the opposite
direction, labelled $\ol{x}$.
Notice that for each path $\pi$ in the loop automaton labelled $w \in \hat{X}^*$,
there is a corresponding path in the opposite direction labelled
$\ol{w}$; we denote this path $\ol{\pi}$.
 We view the loop automaton as a (typically
infinite) automaton over $\hat{X}$, with start state and terminal
state the identity of $M$.
The \textit{(right) loop problem} of $M$ with respect to $\sigma$ is the
language $L_\sigma(M) \subseteq \hat{X}^*$ of words recognised by
the loop automaton $\loopauto{\sigma}{M}$.

Words in $X^*$ and $\ol{X}^*$ are called \textit{positive} and
\textit{negative} words respectively; words in $\hat{X}^*$ which are
neither positive nor negative are called \textit{mixed}. Similarly, an
edge or path in $\loopauto{\sigma}{M}$ is called \textit{positive}
[\textit{negative}, \textit{mixed}] if it has a positive [respectively
negative, mixed] label.

We now introduce the corresponding definitions for more general semigroups,
possibly without identity. Suppose $S$ is a semigroup. By a
\textit{choice of (semigroup) generators} for $S$ we mean a surjective
morphism $\sigma : X^+ \to S$ from a free semigroup onto $S$; again,
$\sigma$ is called \textit{finite} if $X$ is finite.
Let $S^1$ be the semigroup obtained from $S$ by adjoining an identity
element, denoted $1$, even if there already is one. Then $\sigma$
extends uniquely to a monoid choice of generators $\sigma^1 : X^* \to S^1$.
We define the \textit{loop automaton} $\loopauto{\sigma}{S}$ of $S$ with
respect to $\sigma$ to be the loop automaton $\loopauto{\sigma^1}{S^1}$.

Note that if $M$ is a monoid, then it admits choices of generators, and
hence Cayley graphs and loop problems, both as a semigroup and as a monoid.
We shall see below (Proposition~\ref{prop_semimonoid}) that there is a close
relationship between the loop problems of $M$ as a monoid and as a semigroup.

For brevity, we introduce a few notational conventions.
If $X$ is a semigroup [respectively, monoid] generating
\textit{subset} for a semigroup [monoid] $S$, then the inclusion map
$X \to S$ extends naturally to a formal choice of generators
$\sigma : X^+ \to S$
[$\sigma : X^* \to S$].
In this case, we write $\Gamma_X(S)$, $\loopauto{X}{S}$ and
$L_X(S)$ for $\Gamma_\sigma(S)$, $\loopauto{\sigma}{S}$ and
$L_\sigma(S)$ respectively. Moreover, where the choice of generators is
clear, we shall omit the subscript entirely, writing simply $\Gamma(S)$,
$\loopauto{}{S}$ and $L(S)$. Similarly, if only one monoid or semigroup
is under consideration, we simply write $\Gamma_\sigma$, $\hat{\Gamma}_\sigma$
and $L_\sigma$, or even just $\Gamma$, $\hat{\Gamma}$ and $L$.

There is of course a dual notion of the \textit{left Cayley graph} of a
monoid or semigroup with respect to a given generating set, obtained by
considering multiplication by generators on the left. Consequently, one
obtains also the \textit{left loop automaton} (which reads input words from right to
left) and the \textit{left loop problem}. In this paper, we shall restrict
attention to right loop problems. The loss of generality through this
restriction is limited by the fact that the left loop problem of $M$ is
simply the reverse of the right loop problem of the opposite monoid of $M$.

In fact the loop problem is one of two closely related languages which can
naturally be associated to a finitely generated monoid or semigroup.
The other, which we might call the \textit{meeting problem}, is simply
the set of all words of the form $u \ol{v}$ such that $u$ and $v$ are words
over the generators representing the same element of the monoid or semigroup.
The meeting problem is clearly the intersection of the loop problem with the
regular language $X^* \ol{X}^*$; in contrast, the loop problem cannot
be described so easily in terms of the meeting problem. Of the two, then, the
loop problem is the more complex; it is more naturally analogous to the word
problem of a group, and better describes the geometry of the Cayley graph of
the monoid; the meeting problem is typically easier to work with and
more closely related to the word problem in the sense of universal algebra,
but is less geometric and does not directly generalise the usual notion of
a word problem (as a language) for groups. One might reasonably expect that
a language-theoretic restriction on the loop problem will place much stronger
conditions on a semigroup than the same restriction applied to the meeting
problem.

\section{Basic Properties}\label{sec_basic}

In this section we describe some basic properties of the loop problem.
The broad objective is to show that finitely generated monoids and semigroups
which are ``similar'' (that is, which enjoy similar structure) give rise
to languages which are ``similar'' (that is, belong to the same classes
of languages, from amongst the classes most frequently arising in formal
language theory). The technique is to show that when one has two finitely
generated monoids or semigroups which are closely related, one can often
pass between the corresponding languages using simple language-theoretic
operations such as morphisms, inverse morphisms, rational transductions and
Kleene closure.

The following lemma, which follows easily from the definitions, describes the
relationship between paths in the loop automaton $L(M)$ and multiplication
in the monoid $M$.
\begin{lemma}[The Zig Zag Lemma]\label{lemma_zigzag}
Let $\sigma : X^* \to M$ be a choice of generators for a monoid $M$, and
let $x, y \in M$. Let $n \geq 1$
and suppose $u_0, \dots u_{n-1}, v_1, \dots v_n \in X^*$. Then the following
are equivalent:
\begin{itemize}
\item[(i)] $\loopauto{\sigma}{M}$ has a path from $x$ to $y$ labelled $u_0 \ol{v_1} u_1 \ol{v_2} \dots u_{n-1} \ol{v_n}$;
\item[(ii)] there exist $p_0, \dots, p_n \in M$ such that $x = p_0$, $y = p_n$,
and $p_i (u_i \sigma) = p_{i+1} (v_{i+1} \sigma)$ for $0 \leq i < n$.
\end{itemize}
\end{lemma}
\begin{proof}
Suppose first that (i) holds, and let $\pi$ be a path from $x$ to $y$
labelled $u_0 \ol{v_1} u_1 \ol{v_2} \dots u_{n-1} \ol{v_n}$. Let $p_0 = 1$,
and for $1 \leq i \leq n$, define $p_i$ be to be the vertex reached by
while tracing the path $\pi$ from $x$, after reading the prefix
$u_0 \ol{v_1} \dots u_{i-1} \ol{v_i}$. That these vertices have the
required property follows easily from the definitions.

Conversely, if (ii) holds then, since the loop automaton contains the Cayley
graph, there is a path $\pi_i$ from each $p_i$ to $p_i (u_i \sigma)$ labelled
$u_i$ and a path $\phi_{i+1}$ from each $p_{i+1}$ to $p_{i+1} (v_{i+1} \sigma) = p_i (u_i \sigma)$ labelled
$v_{i+1}$. In the loop automaton, each path $\phi_{i+1}$ has an inverse path
$\ol{\phi_{i+1}}$ which runs from $p_i (u_i \sigma) = p_{i+1} (v_{i+1} \sigma)$ to
$p_{i+1}$ and has label $\ol{v_i}$. Now the concatenation
$$\pi_0 \ol{\phi_1} \pi_1 \dots \pi_{n-1} \ol{\phi_{n}}$$
is clearly a path with the properties required to show that (i) holds.
\end{proof}

Since every word in $\ol{X}^*$ can be factored in the form
$u_0 \ol{v_1} u_1 \ol{v_2} \dots u_{n-1} \ol{v_n}$ for some words
$u_0, \dots u_{n-1}, v_1, \dots v_n \in X^*$, the Zig Zag Lemma gives a
complete algebraic description of all the path labels between any
two vertices.
The special case of the Lemma in which $n = 1$ and $x = y = 1$
shows that the loop
problem of a monoid encodes the word problem in the sense of universal
algebra \cite[Section~30]{Graetzer68}. Substituting $S^1$ for $M$ in the
statement gives a corresponding result for semigroups.
\begin{corollary}\label{cor_univalgwp}
Let $\sigma : X^+ \to S$ [respectively, $\sigma : X^* \to M$] be a choice of
generators for a semigroup [monoid] $M$, and let $u, v \in X^+$ [$u,v \in X^*$]. Then
$u \sigma = v \sigma$ if and only if $u \ol{v} \in L_\sigma$.
\end{corollary}
From this, we see immediately that, unlike the sets of words representing
the identity, the loop problem is sufficient to determine 
up to isomorphism the monoid or semigroup, together with its generating system.
\begin{corollary}
Let $X$ be a finite alphabet, and $M_1$ and $M_2$ be monoids with
choices of generators $\sigma_1 : X^* \to M_1$ and $\sigma_2 : X^* \to M_2$.
If $L_{\sigma_1}(M_1) = L_{\sigma_2}(M_2)$ then there is a monoid
isomorphism $\rho : M_1 \to M_2$ such that $\sigma_1 \rho = \sigma_2$.
The corresponding statement for semigroups also holds.
\end{corollary}
\begin{proof}
Suppose $L_{\sigma_1}(M_1) = L_{\sigma_2}(M_2)$. Define $\rho : M_1 \to M_2$ by $(w \sigma_1) \rho = w \sigma_2$.
This map is well-defined since if $u \sigma_1 = v \sigma_1$ then by
Corollary~\ref{cor_univalgwp} we have
$u \ol{v} \in L_{\sigma_1}(M_1)$; but then $u \ol{v} \in L_{\sigma_2}(M_2)$
so by Corollary~\ref{cor_univalgwp} again we have $u \sigma_2 = v \sigma_2$.
An entirely similar argument establishes that $\rho$ is injective.

To see that $\rho$ is a monoid homomorphism, we observe that
$$[(u \sigma_1) (v \sigma_1)] \rho = ((uv) \sigma_1) \rho = (uv) \sigma_2
= (u \sigma_2) (v \sigma_2) = [(u \sigma_1) \rho] [(v \sigma_1) \rho].$$
Surjectivity of $\rho$ is immediate from surjectivity of $\sigma_2$.
Finally, that $\sigma_1 \rho = \sigma_2$ follows straight from the
definition, completing the result for monoids.

The statement for semigroups follows directly from the definition and the
result for monoids, with $S^1$ substituted for $M$.
\end{proof}

\begin{lemma}\label{lemma_loopsdown}
Let $\sigma : X^* \to M$ be a choice of generators for a monoid $M$, and
suppose $w \in L_\sigma(M)$. Then $w$ labels a loop at every vertex in
$\loopauto{\sigma}{M}$. The corresponding statement for semigroups also
holds.
\end{lemma}
\begin{proof}
Let $x \in M$.
Write $w = u_0 \ol{v_1} u_1 \ol{v_2} \dots u_{n-1} \ol{v_n}$ with each
$u_i, v_i \in X^*$. That $w \in L_\sigma{M}$ means exactly that
$w \in \hat{X}^*$ labels a loop at $1$ in $\loopauto{\sigma}{M}$.
By the Zig Zag Lemma (Lemma~\ref{lemma_zigzag}),
there exist $p_0, \dots, p_n \in M$ such that $p_0 = p_n = 1$ and
$p_i (u_i \sigma) = p_{i+1} (v_{i+1} \sigma)$ for
$0 \leq i < n$. For each $i$, let $q_i = x p_i$. Then we have $p_0 = p_1 = x$,
and for each $i$, $q_i (u_i \sigma) = x p_i (u_i \sigma) = x p_{i+1} (v_{i+1} \sigma) = q_{i+1} (v_{i+1} \sigma)$. 
Now applying the Zig Zag Lemma again, we obtain the required loop at
$x$ labelled $w$.

The corresponding result for semigroups is once again obtained by
substituting $S^1$ for $M$, and applying the definition of the loop
automaton of a semigroup.
\end{proof}

We apply Lemma~\ref{lemma_loopsdown} to prove a simple language-theoretic
property of the loop problem. Recall that a language $L$ is
\textit{insertion-closed} if whenever $x, y, w$ are (possibly empty) words
such that $xy \in L$ and
$w \in L$ we have $xwy \in L$.
\begin{proposition}\label{prop_insertclosed}
Every loop problem is insertion-closed.
\end{proposition}
\begin{proof}
By definition, every loop problem of a semigroup $S$ is also a loop problem of
the monoid $S^1$ with respect to a choice of monoid generators, so
it suffices to consider the loop problem $L_\sigma(M)$ of a monoid $M$ with
respect to some choice of monoid generators $\sigma : X^* \to M$.
Suppose $xy \in L_\sigma(M)$ and $w \in L_\sigma(M)$. Then $\loopauto{\sigma}{M}$ has a vertex $v$ such
that there is a path from $1$ to $v$ labelled $x$ and a path from $v$ to $1$ labelled $y$. Now by
Lemma~\ref{lemma_loopsdown}, there is a loop at $v$ labelled $w$. Hence,
there is a path from $1$ to $1$ labelled $xwy$ so $xwy \in L_\sigma(M)$.
Thus, $L_\sigma(M)$ is insertion-closed, as required.
\end{proof}

It is well-known that the word problems of a group with respect to different
finite choices of generators can be obtained from each other as inverse
images under morphisms of free semigroups. Many important classes of
languages are closed under inverse morphism, and it follows that these
classes correspond to invariant properties of groups which do
not depend upon the finite choice of generators. It transpires that the
loop problems of monoids and semigroups enjoy a similar property.
\begin{proposition}\label{prop_gens}
Let $M$ be a monoid and $\sigma : X^* \to M$ and $\tau : Y^* \to M$
be finite choices of generators. Then $L_\sigma(M)$ is an inverse morphic
image of $L_\tau(M)$. The corresponding statement for semigroups also holds.
\end{proposition}
\begin{proof}
For each generator $y \in Y$, let $w_y$ be a word in $X^*$ such that
$w_y \sigma = x \tau$. Define a morphism $\rho : \hat{Y}^* \to \hat{X}^*$
by $y \rho = w_y$ and $\ol{y} \rho = \ol{w_y}$ for all $y \in Y$. Note that
$u \rho \sigma = u \tau$ for every $u \in Y^*$, and $\ol{v} \rho = \ol{v \rho}$
for every $v \in \hat{Y}^*$. We claim that
$L_\tau(M) = L_\sigma(M) \rho^{-1}$.

Suppose $w \in \hat{Y}^*$. Write
$$w = u_0 \ol{v_1} u_1 \ol{v_2} \dots u_{n-1} \ol{v_n}$$
where
$u_0, \dots u_{n-1}, v_0, \dots, v_n \in Y^*$. Then $w \in L_\tau(M)$ if and
only if the loop automaton $\loopauto{\tau}{M}$
has a path from $1$ to $1$ labelled $w$. By the Zig Zag Lemma (Lemma~\ref{lemma_zigzag}), this is
true exactly if there exist elements
$p_0, \dots, p_n \in M$ such that $p_0 = p_n = 1$
and $p_i (u_i \tau) = p_{i+1} (v_{i+1} \tau)$ for $0 \leq i < n$.
But this holds exactly if there exist elements
$p_0, \dots, p_n \in M$ such that $p_0 = p_n = 1$
and $p_i (u_i \rho \sigma) = p_{i+1} (v_{i+1} \rho \sigma)$ for
$0 \leq i < n$. Now by the Zig Zag Lemma again,
this is true if and only if
\begin{align*}
w \rho &= (u_0 \rho) (\ol{v_1} \rho) (u_1 \rho) (\ol{v_2} \rho) \dots (u_{n-1} \rho) (\ol{v_n} \rho) \\
&= (u_0 \rho) (\ol{v_1 \rho}) (u_1 \rho) (\ol{v_2 \rho}) \dots (u_{n-1} \rho) (\ol{v_n \rho})
\end{align*}
labels a loop at $1$ in $\loopauto{\sigma}{M}$, that is, lies in $L_\sigma(M)$.

For the semigroup case, if $\sigma : X^+ \to S$ and $\tau : Y^+ \to S$ are
semigroup choices of generators, then $\sigma^1 : X^+ \to S^1$ and
$\tau^1 : Y^* \to S^1$ are monoid choices of generators for the monoid
$S^1$, so arguing as above, we obtain a morphism $\rho : Y^* \to X^*$ such
that
$$L_\tau(S) = L_{\tau^1}(S^1) = L_{\sigma^1}(S^1) \rho^{-1} = L_\sigma(S)$$
as required.
\end{proof}

Proposition~\ref{prop_gens} tells us that, for example, the classes of
regular and context-free languages can be expected to correspond to
invariant properties of finitely generated monoids. It is not, however,
immediately obvious that these properties will be natural ones from a
semigroup-theoretic perspective. In Section~\ref{sec_correspond} below,
we shall attempt to isolate the properties in question.

Suppose $M$ is a monoid. We remarked in Section~\ref{sec_basic} that
$M$ admits loop problems both as a semigroup and a monoid.
Proposition~\ref{prop_gens} ensures that the different loop problems for
$M$ as a monoid are closely related, as are those for $M$ as a semigroup.
It does not, however, guarantee a relationship between the loop problems
with respect to monoid generating sets, and those with respect to semigroup
generating sets. The following proposition shows that there is nonetheless
a close relationship.
\begin{proposition}\label{prop_semimonoid}
Let $M$ be a monoid, let $\sigma : X^+ \to M$ be a semigroup choice
of generators for $M$, and let $\tau : X^* \to M$ be the unique extension
of $\sigma$ to a monoid choice of generators for $M$. Then
$$L_\sigma(M) = L_\tau(M) \cap (X \hat{X}^* \ol{X} \cup \lbrace \epsilon \rbrace),$$
and there exists a word $w \in X^+$ such that
$$L_\tau(M) = w^{-1} L_\sigma(M) \ol{w}^{-1} = \lbrace u \in \hat{X}^* \mid w u \ol{w} \in L_\sigma(M) \rbrace.$$
\end{proposition}
\begin{proof}
Since by definition $L_\sigma(M) = L_{\sigma^1}(M^1)$, we need to consider
both the identity of $M$ and the extra identity adjoined in $M^1$; we write
$e$ for the identity of $M$, and $1$ for the identity of $M^1$. Observe
that $\loopauto{\sigma}{M} = \loopauto{\sigma^1}{M^1}$ is isomorphic to
$\loopauto{\tau}{M}$ with an extra vertex $1$ adjoined. The edges incident
with $1$ are as follows; for each generator $x \in X$, there is an edge from $1$ to $x \sigma$ labelled $x$, and of
course a corresponding inverse edge from $x \sigma$ to $1$ labelled $\ol{x}$.
In particular, whenever there is an edge from $1$ to a vertex $v$ [respectively,
$v$ to $1$], there is also an edge in both $\loopauto{\sigma}{M}$ and
$\loopauto{\tau}{M}$ from $e$ to $v$ [$v$ to $e$] with the same label.

We prove first that
$$L_\sigma(M) = L_\tau(M) \cap (X \hat{X}^* \ol{X} \cup \lbrace \epsilon \rbrace).$$
It is immediate from our description above of the edges incident with $1$,
that every non-empty label of a loop at $1$ in $\loopauto{\sigma}{M}$ must begin with a positive
generator and end with a negative generator, so that $L_\sigma(M) \subseteq X \hat{X}^* \ol{X} \cup \lbrace \epsilon \rbrace.$
Moreover, from our observation that every edge incident with $1$ corresponds
to an edge incident with $e$ in $\loopauto{\tau}{M}$, it follows
easily that every word accepted by $\loopauto{\sigma}{M}$ is accepted by
$\loopauto{\tau}{M}$, so that $L_\sigma(M) \subseteq L_\tau(M)$.

Conversely, suppose $u \in L_\tau(M) \cap (X \hat{X}^* \ol{X} \cup \lbrace \epsilon \rbrace).$
If $u = \epsilon$ then certainly $u \in L_\sigma(M)$, so we may assume that
$u \in X \hat{X}^* \ol{X}$, that is, $u = x v \ol{y}$ for some $x,y \in X$
and $v \in \hat{X}^*$. Let $\pi$ be a loop at $e$ labelled $u$ in the
loop automaton $\loopauto{\tau}{M}$.
 Then $\pi$ begins with an edge leaving $e$ with label $x$;
from the definition of the loop automaton, it follows easily that
the target of this edge is the vertex $x \sigma$. Similarly, $\pi$ ends
with an edge from $y \sigma$ to $e$ labelled $\ol{y}$. In between is a
path from $x \sigma$ to $y \sigma$ labelled $v$.
But $\loopauto{\sigma}{M}$ has an edge from $1$ to $x \sigma$ labelled $x$,
and an edge from $y \sigma$ to $1$ labelled $\ol{y}$. Moreover, since
$\loopauto{\sigma}{M}$ contains $\loopauto{\tau}{M}$, it also has a path
from $x \sigma$ to $y \sigma$ labelled $v$. Thus, we obtain a path in $\loopauto{\sigma}{M}$ from $1$ to $1$
with label $u$, so that $u \in L_\sigma(M)$. This completes the proof of the
first claim.

Now let $w \in X^+$ be any non-empty word such that $w \sigma = e$. We
claim that 
$$L_\tau(M) = w^{-1} L_\sigma(M) \ol{w}^{-1} = \lbrace u \in \hat{X}^* \mid w u \ol{w} \in L_\sigma(M) \rbrace.$$

To show one inclusion, suppose $u \in L_\tau(M)$, so that $u \in \hat{X}^*$ labels a loop at $e$ in
$\loopauto{\tau}{M}$.
Clearly $\loopauto{\sigma}{M}$ has a path from $1$ to $e$ labelled $w$ and
hence a path from $e$ to $1$ labelled $\ol{w}$. Since 
$\loopauto{\sigma}{M}$ contains $\loopauto{\tau}{M}$, it also has a loop
at $e$ labelled $u$. Thus, $\loopauto{\sigma}{M}$ has a
loop at $1$ labelled $w u \ol{w}$, so that
$u \in w^{-1} L_\sigma(M) \ol{w}^{-1}$.

Conversely, suppose that $u \in w^{-1} L_\sigma(M) \ol{w}^{-1}$, that is, that
$wu\ol{w} \in L_\sigma(M)$, so that $\loopauto{\sigma}{M}$ has a loop at $1$
labelled $w u \ol{w}$. It follows easily from the definition of the loop
automaton that there is a unique path beginning at $1$ with label $w$, and that
this path ends at $e$. Dually, there is a unique path ending at $1$ with
label $\ol{w}$, and this path begins at $e$.
Hence, $\loopauto{\sigma}{M}$ must have a loop at
$e$ labelled $u$. Moreover, because of our observation at the start of the
proof that every edge incident with $1$ has a corresponding edge incident
with $e$, we can find such a loop which does not visit the vertex $1$.
It follows that $\loopauto{\tau}{M}$ also has a loop at $e$ labelled $u$,
so that $u \in L_\tau(M)$, as required to complete the proof of the second
claim.
\end{proof} 

Since taking an intersection with a regular languages, left and right
translation, and taking an inverse image under a morphism are all operations
which can be performed by rational transductions \cite{Berstel79}, 
Propositions~\ref{prop_gens} and~\ref{prop_semimonoid} yield the following
result, which says that the loop problem of a semigroup or monoid is, modulo
rational equivalence, invariant under choice of generators.
\begin{corollary}\label{cor_gens}
Let $S$ be a semigroup or monoid. Then any two loop problems for $S$ (as
a semigroup or as a monoid if appropriate) are rationally equivalent.
\end{corollary}
In view of Corollary~\ref{cor_gens}, we  shall permit ourselves to speak
simply of ``the'' loop problem of a monoid  or semigroup where, as is often the case, our interest is in 
language-theoretic properties which are invariant under rational 
transductions.

\section{Groups and Completely Simple Semigroups}\label{sec_groupcs}

In this section, we consider the loop problems of groups, and then more
generally of completely simple semigroups.

Let $\sigma : X^* \to G$ be a choice of monoid generators for a group $G$.
Recall that the \textit{word problem} for $G$ with respect to $\sigma$ is
the language of all words $w \in X^*$ such that $w \sigma = 1$ in $G$.
Equivalently, it is the language accepted by the Cayley graph $\Gamma_\sigma(G)$
when viewed as an automaton with initial and terminal state the identity of
$G$. An obvious question concerns
the relationship between the loop problem for a group (with respect to a
particular finite choice of generators) and the word problem. The
following proposition says that they are almost the same thing.
\begin{proposition}\label{prop_grouploop}
Let $\sigma : X^* \to G$ be a choice of monoid generators for a group $G$.
Then there exists a choice of monoid generators $\tau : \hat{X}^* \to G$ such that the loop
problem $L_\sigma(G)$ of $G$ with respect to $\sigma$ is equal to the word
problem of $G$ with respect to $\tau$.
\end{proposition}
\begin{proof}
Define $\tau : \hat{X}^* \to G$ by
$x \tau = x \sigma$ and $\ol{x} \tau = (x \sigma)^{-1}$ for all
$x \in X$. It is readily verified that $\loopauto{\sigma}{G}$ is exactly the same as
$\Gamma_\tau(G)$, so that $L_\sigma(G)$ is equal to the word problem of $G$
with respect to $\tau$, as required.
\end{proof}
Since we are chiefly interested in intrinsic structural properties of
groups (and semigroup and monoids) which are not dependent upon a particular
choice of generators, this correspondence is quite sufficient for our
purposes.

Next, we consider the broader class of \textit{completely simple} semigroups.
Recall that a \textit{primitive idempotent} in a semigroup $S$ is an idempotent $e$
with the property that for any non-zero idempotent $f$
such that $ef = fe = f$, we have $e = f$. A semigroup is
called \textit{completely simple} if it has a primitive idempotent and no
proper ideals. For a detailed introduction to the theory of completely
simple semigroups, including a number of equivalent definitions, see
\cite[Chapter~2]{Clifford61} or \cite[Chapter~3]{Howie95}. A
\textit{subgroup} of a semigroup $S$
is a subsemigroup of $S$ which forms a group with the multiplication
inherited from $S$; a subgroup of $S$ is
\textit{maximal} if it is not properly contained in any other subgroup.
A construction of Rees \cite{Rees40} gives a
simple combinatorial description of completely simple semigroups, in
terms of their maximal subgroups. Let $G$ be a group and $I$ and
$J$ be sets, and let $P$ be a $J \times I$ matrix with entries
drawn from $G$. The \textit{Rees matrix semigroup} $M(G; I, J; P)$ is the
semigroup with set of elements $I \times G \times J$ and multiplication
given by $(i, g, j) (k, h, l) = (i, g P_{jk} h, l)$ for all $i,k \in I$,
$g,h \in G$ and $j,l \in J$. The following theorem is usually attributed
to Rees, although the essential idea was given by Suschkewitz \cite{Suschkewitz28}.
\begin{theorem}[Suschkewitz 1928, Rees 1940]\label{thm_rees}
Let $G$ be a group, $I$ and $J$ sets, and $P$ a
$J \times I$ matrix over $G$. Then the Rees matrix semigroup
$M(G; I, J; P)$ is completely simple with maximal subgroups all isomorphic
to $G$. Conversely, every completely simple semigroup is isomorphic to one
of this form.
\end{theorem}
Both completely simple semigroups and Rees matrix constructions are of
central importance in the structural theory of semigroups. We shall need
the following elementary property of completely simple semigroups.

\begin{lemma}\label{lem_csprops}
Let $H$ be a maximal subgroup of a completely simple semigroup $S$. 
If $a \in S^1$ and $x,y \in H$ are such that $ax \in H$, then
there exists $b \in H$ such that $bx = ax$ and $by = ay$.
\end{lemma}
\begin{proof}
To prove the claim when $a = 1$ we simply take $b$ to be the identity of $H$.
Otherwise, by the Rees theorem as described above, we may assume that
$S = M(G; I, J; P)$ for some group $G$, sets $I$ and $J$ and $J \times I$
sandwich matrix $P$. It is well known, and can easily be deduced from the
Rees theorem, that there exist $i \in I$ and $j \in J$ such that
the maximal subgroup $H$ is the set of all elements of the form $(i, g, j)$ for
different $g \in G$.

Thus, we may suppose that $a = (i_a, g_a, j_a)$, $x = (i, g_x, j)$ and
$y = (i, g_y, j)$.
Set $b = (i, g_a P_{j_a i} P_{ji}^{-1}, j) \in H$. Now we have $b \in H$ and
$$ax = (i_a, g_a, j_a) (i, g_x, j) = (i_a, g_a P_{j_a i} g_x, j).$$
On the other hand,
$$bx = (i, g_a P_{j_a i} P_{ji}^{-1}, j) (i, g_x, j) = (i, g_a P_{j_a i} P_{ji}^{-1} P_{ji} g_x, j) = (i, g_a P_{j_a i} g_x, j).$$
Moreover, since we know that $ax \in H$, we must have $i_a = i$, from which
it follows that $ax = bx$ as required. An similar argument (but using the fact
that we already know that $i_a = i$ in order to avoid the need to presuppose
that $ay \in H$) shows that $ay = by$, completing the proof.
\end{proof}

The Rees theorem often allows results about
groups to be extended in some form to completely simple semigroups. The
following two theorems say that the loop problem of a finitely generated
completely simple semigroup is closely related to the loop problem (or
equivalently the word problem) of its maximal subgroups.

\begin{theorem}\label{thm_cswholetosub}
Let $G$ be a maximal subgroup of a completely simple semigroup $S$.
Let $\sigma : X^+ \to G$ be a choice of semigroup generators for $G$, and
$\tau : Y^+ \to S$ a choice of generators for $S$, such that $X \subseteq Y$
and $\sigma$ is the restriction of $\tau$ to $X^+$. Then
$L_\sigma(G) = L_\tau(S) \cap \hat{X}^*$.
\end{theorem}
\begin{proof}
The
loop automaton $\loopauto{\sigma}{G}$ of $G$ is clearly embedded in a natural
way into the loop automaton $\loopauto{\tau}{S}$. It follows easily
that any word accepted by the former is also accepted by the latter, so
that $L_\sigma(G) \subseteq L_\tau(S) \cap \hat{X}^*$.

Conversely, suppose $w$ is a word in $L_\tau(S) \cap \hat{X}^*$. Then
$\loopauto{\tau}{S} = \loopauto{\tau^1}{S^1}$ has a loop at $1$ labelled
$w$. Write $w = u_0 \ol{v_1} u_1 \ol{v_2} \dots u_{n-1} \ol{v_n}$ with each
$u_i, v_i \in X^*$. By the Zig Zag Lemma (Lemma~\ref{lemma_zigzag}),
there exist $p_0, \dots, p_n \in S^1$ such that $p_0 = p_n = 1$ and
$p_i (u_i \tau) = p_{i+1} (v_{i+1} \tau)$ for
$0 \leq i < n$. We claim that the elements $p_0, \dots, p_{n-1}$ can all be
chosen to lie in $G^1$. Indeed, suppose we are given $p_0, \dots, p_{n-1}$
satisfying the above equations and not all lying in $G^1$, and let $k$ be
minimal such that $p_k$ does not lie in $G^1$. Certainly, 
$1 \leq k \leq n-1$.
Now we have $u_{k-1} \tau \in G$ and $p_{k-1} \in G^1$, so that
$p_{k-1} (u_{k-1} \tau) = p_k (v_k \tau)$ lies in $G$. Since $v_k \tau$ and
$u_k \tau$ also lie in $G$, it follows by Lemma~\ref{lem_csprops} that
there exists an element $q_k \in G$ with
$q_k (v_k \tau) = p_k (v_k \tau) = p_{k-1} (u_{k-1} \tau)$ and
$q_k (u_k \tau) = p_k (u_k \tau) = p_{k+1} (v_{k+1} \tau)$. Hence, we
can replace $p_k$ with $q_k$ to obtain a sequence $p_0, \dots, p_n$ with
strictly fewer elements outside $G^1$. Continuing in the same vein, we
eventually obtain $p_0, \dots, p_n \in G^1$ with the desired properties.

Applying the Zig Zag Lemma again, we now see that
$\loopauto{\sigma}{G} = \loopauto{\sigma^1}{G^1}$ has a
loop at $1$ labelled $w$, so that $w \in L_\tau(G)$, as required.
\end{proof}

\begin{theorem}\label{thm_cssubtowhole}
Let $S$ be a finitely generated completely simple semigroup with maximal
subgroups isomorphic to a group $G$. Then the loop problem for $S$ is the
Kleene closure of a rational transduction of the word problem for $G$.
\end{theorem}
\begin{proof}
By the Rees theorem, as described above, we may assume that
$S = M(G; I, J; P)$ where by the main theorem of \cite{Ayik99}, $I$ and
$J$ are finite.
Let $\sigma : X^* \to G$ be a finite choice of monoid generators for $G$,
and $\tau : Y^+ \to M$ a finite choice of semigroup generators for $S$.
For each $y \in Y$, suppose $y \tau = (i_y, g_y, j_y)$ and let
$w_y, w_y' \in X^*$
be words representing $g_y, g_y^{-1} \in G$ respectively. For each $i \in I$
and $j \in J$, let $w_{ji}, w_{ji}' \in X^*$ be words representing
$P_{ji}, P_{ji}^{-1} \in G$ respectively.

We define a finite state transducer from $X$ to $Y$ with
\begin{itemize}
\item vertex set $(I \times J) \ \cup \ \lbrace A, Z \rbrace$ where $A$ and $Z$
are new symbols;
\item initial state $A$;
\item terminal state $Z$;
\item for each generator $y \in Y$, an edge from $A$ to
$(i_y,j_y)$ labelled $(w_y, y)$;
\item for each generator $y \in Y$, an edge from $(i_y,j_y)$ to
$Z$ labelled $(w_y', \ol{y})$;
\item for each generator $y \in Y$ and each $k \in J$, an
edge from $(i_y,k)$ to $(i_y,j_y)$ labelled $(w_{ki_y} w_y, y)$; and
\item for each generator $y \in Y$ and each
$k \in J$, an edge from $(i_y,j_y)$ to $(i_y,k)$ labelled
 labelled $(w_y' w_{ki_y}', \ol{y})$.
\end{itemize}

Now let $g \in G$, $i \in I$, $j \in J$ and $v \in \hat{Y}^*$ and suppose
$n$ is a positive integer. We say that a path in the loop automaton
$\loopauto{\tau}{S}$ which originates at $1$ is \textit{non-returning} if
it does not visit the vertex $1$ at any point other than the start and
(possibly) the end. As a first step towards the proof, we claim that
the following conditions are equivalent.
\begin{itemize}
\item[(i)] the loop automaton $\loopauto{\tau}{S}$ has a non-returning path
           of length $n$ from $1$ to $(i, g, j)$ labelled $v$;

\item[(ii)] the transducer has a path of length $n$ from $A$ to $(i,j)$
            labelled $(u, v)$ for some $u \in X^*$ which represents $g$.
\end{itemize}
The proof of equivalence proceeds by induction on the path length $n$. That
the equivalence holds when $n = 1$ follows immediately from the definition
of the edges in the transducer. Now let $n > 1$ and suppose true for smaller $n$.

Suppose first that (i) holds, and let $\pi$ be the path given by the
hypothesis. Let $e$ be the last edge the path $\pi$, and let $\pi'$ be
the path $\pi$ with the last edge removed, so that $\pi = \pi' e$. Let
$v'$ be the label of $\pi'$. Since $n > 1$ and $\pi$ is non-returning,
the path $\pi'$ must end at a vertex
of the form $(i',g',k)$.
It follows easily from the definition of
the multiplication in a Rees matrix semigroup that the vertices in the
loop automaton corresponding to elements with first coordinate $i$ are
connected to the rest of the automaton only via the vertex $1$.
Hence, since the path $\pi$ is non-returning, we must have $i = i'$, so
that $\pi$ actually ends at $(i,g',k)$. Now $\pi$ is a path of length $n-1$, so by
the inductive hypothesis, the transducer has a path of length $n-1$ from
$A$ to $(i,k)$ labelled $(u', v')$ for some word $u' \in X^*$ which
represents $g' \in G$.

We now treat separately the case where $e$ is a positive edge, and that
where $e$ is a negative edge. Suppose first that $e$ is a positive edge,
with label $y \in Y$. Then from the definition, the transducer
has an edge from $(i, k)$ to $(i,j)$ with label $(w_{ki_y} w_y, y)$. Hence,
the transducer has a path of length $n$ from $1$ to $(i,j)$ with label
$$(u', v') (w_{ki_y} w_y, y) \ = \ (u' w_{ki_y} w_y, v' y) \ = \ (u' w_{ki_y} w_y, v).$$
Now from the definition of the loop automaton, we must have
$$(i, g, j) \ = \ (i,g',k) (y \tau) \ = \ (i, g', k) (i_y, g_y, j_y) \ = \ (i, g' P_{k i_y} g_y, j_y).$$
Equating second coordinates, we see that $g = g' P_{ki_y} g_y$. But it
follows that $g$ is represented by the word $u' w_{ki_y} w_y \in X^*$, so
setting $u = w_{ki_y} w_y \in X^*$, we see that (ii) holds as required.

On the other hand, suppose $e$ is a negative edge, with label $\ol{y}$ for
some $y \in Y$. In this case the transducer by definition has an edge from
$(i,k)$ to $(i,j)$ with label
$(w_y' w_{ki_y}', \ol{y})$. Hence, there is a path of length $n$
from $1$ to $(i,j)$ with label
$$(u', v') (w_y' w_{ki_y}', \ol{y}) \ = \ (u' w_y' w_{ki_y}', v' \ol{y}) \ = \ (u' w_y^{-1} w_{ki_y}^{-1}, v).$$
Now from the definition of the loop automaton, we must have
$$(i, g', j) \ = \ (i,g',k) (y \tau) \ = \ (i, g', k) (i_y, g_y, j_y) \ = \ (i, g' P_{k i_y} g_y, j_y).$$
Again equating second coordinates, we see this time that $g P_{ki_y} g_y = g'$, so
that $g = g' g_y^{-1} P_{ki_y}^{-1}$, and $g$ is represented by the word
$u' w_y' w_{ki_y}'$. Now setting $u = u' w_y^{-1} w_{ki_y}^{-1}$,
we again see that (ii) holds as required. This completes the proof that
(i) implies (ii).

Conversely, suppose (ii) holds, and this time let $\pi$ be a path of length
$n$ in the transducer from $A$ to $(i,j)$ labelled $(u,v)$ for some $u \in X^*$
which represents $g$. Much as before, we let $e$ be the last edge of $\pi$ and
$\pi'$ be the path $\pi$ with the final edge removed. Then $\pi'$ is
a path of length $n-1$ from $A$ to some vertex $(i',k)$ with label of the
form $(u', v')$. Moreover, it follows easily from the definition of the
transducer that $i = i'$, so that $\pi'$ ends at $(i, k)$.
Let $g' \in G$ be the element represented by $u'$. Then
by the inductive hypothesis,
there exists a path of length $n-1$ in the loop automaton from $1$ to
$(i,g',k)$ with label $v'$. Now $e$ is an edge from $(i,k)$ to
$(i,j)$. From the definition of the edges in the transducer, we see that there
exists $y \in Y$ with $j_y = j$ such that $e$ has label either
$(w_{ki} w_y, y)$ or $(w_y' w_{ki_y}', \ol{y})$. As before, we treat these
two cases
separately. In the former case, observe that we have $u = u' w_{ki_y} w_y$
from which we deduce that $g = g' P_{ki_y} g_y$. But now
$$(i,g',k) (y \sigma) = (i, g' P_{ij_y} g_y, j_y) = (i, g, j)$$
so we see that the loop automaton has an edge from $(i,g',k)$ to
$(i,g,j)$ labelled $y$. Combining this with the path whose existence we deduced
using
the inductive hypothesis, we conclude that the loop automaton has a path
from $1$ to $(i,g,j)$ labelled $v = v' \ol{y}$, so that (i) holds as required.
An entirely similar argument suffices to show that (i) also holds in the
case that $e$ has label of the form  $(w_y' w_{ki_y}', \ol{y})$, thus
completing the proof of the equivalence of conditions (i) and (ii)
above.

Now let $K$ be the language of all words which label non-returning loops at
the identity in $\loopauto{Y}{S}$. We claim now that $K$ is exactly the
image of the word
problem of $G$ under the transduction defined by our transducer.
Clearly the Kleene closure $K^*$ is exactly the loop problem $L_\tau(S)$,
so this will suffice to complete the proof of the theorem.

Suppose first that $v \in K$. Then by definition the loop automaton
$\loopauto{Y}{S}$ has a non-returning loop at $1$ labelled $v$. Note that
all edges in $\loopauto{Y}{S}$ which end at $1$ run from vertices
corresponding to generators
$y$ and have label $\ol{y}$, so we may assume that the last edge of the path
runs from a generator $y$ to $1$, and has label $\ol{y}$. Let $\pi$ be the path
without this last edge, so that $\pi$ runs from $1$ to $y$, and
let $v' \in \hat{Y}^*$ be the label of this path, so that $v = v' \ol{y}$.
Then by the equivalence
above, the transducer has a path from $A$ to $(i_y, j_y)$ labelled $(u,v')$
for some $u \in X^*$ representing $g_y$. But directly from the definition,
the transducer also has an edge from $(i_y, j_y)$ to $Z$ labelled
$(w_y^{-1}, \ol{y})$.
Hence, we deduce that $(u w_y^{-1}, v)$ is accepted by the transducer,
where $u w_y^{-1}$ represents $g_y g_y^{-1} = 1$ in $G$, so that $v$
lies in the image under the transduction of the word problem of $G$.

Conversely, suppose $u$ is a word representing $1$ in $G$, such that the
transducer accepts $(u,v)$. Then the transducer has a path $\pi$ from $A$
to $Z$ labelled $(u,v)$. Again, we proceed by letting $\pi'$ be the path
obtained from $\pi$ by deleting the last edge. Then there must exist a
generator $y$ such that $\pi'$ ends at $(i_y,j_y)$. Moreover, $\pi'$ must
be labelled $(u', v')$ where $u = u' w_y'$ and $v = v' \ol{y}$. Now
since $u$ represents $1$ and $w_y'$ represents $g_y^{-1}$, we deduce that
$u'$ represents $g_y$. By the equivalence above, it follows that the
loop automaton $\loopauto{\tau}{S}$ has a non-returning path from $1$ to
$(i_y, g_y, j_y)$ labelled $v'$. Now it certainly also
has an edge from $(i_y, g_y, j_y)$ to $1$ labelled $\ol{y}$, so we deduce
that $v = v' \ol{y} \in K$, as required. This completes the proof.
\end{proof}

We remark that the transducer constructed in the proof of 
Theorem~\ref{thm_cssubtowhole} can be construed as a \textit{$G$-automaton}
\cite{KambitesGAuto} accepting part of the loop problem for $S$.

Combining Theorems~\ref{thm_cswholetosub} and \ref{thm_cssubtowhole} with
Propositions~\ref{prop_semimonoid} and \ref{prop_grouploop}, we obtain the following.
\begin{theorem}\label{thm_csmaxsubgp}
Let $\FamilyF$ be a family of languages closed under rational transduction and
Kleene closure, and let $S$ be a finitely generated completely simple
semigroup. Then the following are equivalent
\begin{itemize}
\item[(i)] the loop problem for $S$ belongs to $\FamilyF$;
\item[(ii)] the loop problem for each maximal subgroup of $S$ belongs to $\FamilyF$;
\item[(iii)] the word problem for each maximal subgroup of $S$ belongs to $\FamilyF$.
\end{itemize}
\end{theorem}

\section{Structural-Linguistic Correspondences}\label{sec_correspond}

In this section, we exhibit some correspondences between the structural
properties of a monoid, and the linguistic properties of its loop problem.
As a first step, we obtain a generalisation of a foundational result of
Anisimov \cite{Anisimov71} from combinatorial group theory.
\begin{theorem}
Let $\sigma : X^* \to M$ be a finite choice of generators for a monoid
$M$. Then $M$ is finite if and only if $L_\sigma(M)$ is regular. The
corresponding statement for semigroups also holds.
\end{theorem}
\begin{proof}
If $M$ is finite then it follows immediately from the definitions that
the loop automaton $\loopauto{\sigma}{M}$ is finite. Hence, $L_\sigma(M)$
is recognised by a finite automaton, and so is regular.

For the converse, recall that a \textit{(right) cone} of a language
$L \subseteq \ol{X}^*$
is a language of the form
$$w^{-1} L = \lbrace x \in \ol{X}^* \mid wx \in L \rbrace$$
for some $w \in \ol{X}^*$. The Myhill-Nerode Theorem \cite[Theorem~3.9]{Hopcroft69}
states that a language is regular if and only if it has finitely many distinct
right cones. Now let $u, v \in X^*$ be words in the generators of $M$. Then
by Corollary~\ref{cor_univalgwp}, we have $u\ol{v} \in L$ if and only if $u$
and $v$ represent the same element of $M$. It follows that $\ol{v} \in u^{-1} L(M)$ if and only if $u$ and
$v$ represent the same element of $M$, and hence that $L(M)$ has a
distinct cone for each element of $M$. So if $M$ is infinite then $L(M)$
has infinitely many cones, and so is not regular.

The corresponding statement for semigroups follows from the result for
monoids together with the fact that a semigroup $S$ is finite if and only
if $S^1$ is finite.
\end{proof}

We now turn our attention to the class of monoids whose loop problem is
a \textit{context-free} language (see \cite{Berstel79}). In the group case,
a well-known theorem of Muller and
Schupp \cite{Muller83}, augmented by a subsequent result of Dunwoody \cite{Dunwoody85},
says that a finitely generated group has context-free word problem if and
only if it is \textit{virtually free}, that is, has a free subgroup of
finite index. Despite the straightforward combinatorial nature of this
statement, the only known proof depends essentially on deep geometric
results about group Cayley graphs, and specifically on Stallings' theory
of \textit{ends} \cite{Stallings68}. We pose the following question.

\begin{question}
Can one characterize or even classify the monoids and semigroups with
context-free
loop problem?
\end{question}
We suspect this question (in general) to be difficult, and believe that
a satisfactory answer and the techniques used to obtain it may be a
significant development with implications reaching well beyond combinatorial
semigroup theory. Firstly, such a result may lead to a purely combinatorial
understanding of the Muller-Schupp Theorem, which has been much sought after
by combinatorial group theorists. Secondly, it may require the development
of new tools for showing that certain languages are \textit{not} context-free;
since there are a number of notable languages conjectured, but not proven,
to be non-context-free, such techniques are likely to be of considerable
interest to formal language theorists

We hope that seeking answers in different classes of ``well-behaved'' 
semigroups may give an indication of how one could proceed in the general 
case. In semigroup theory, of course, ``well-behaved'' typically means one 
of ``combinatorially straightforward'', ``in some way group-like'' or
perhaps ``combinatorially straightforward modulo 
group theory''. One might expect that the question can be answered in the 
first case by elementary combinatorial means, and in the last by application
of the Muller-Schupp theorem. As one example, recalling that the class of
context-free languages is closed under rational transduction \cite[Corollary~4.2]{Berstel79} and 
Kleene closure \cite[Theorem~2.1]{Berstel79}, we can combine the 
Muller-Schupp Theorem with Theorem~\ref{thm_csmaxsubgp} to obtain a 
complete description of completely simple semigroups with context-free 
loop problem.
\begin{theorem}
Let $S$ be a finitely generated completely  simple semigroup. Then $S$ has
context-free loop problem if and only if its maximal subgroups are virtually
free.
\end{theorem}

To conclude this section, we discuss briefly the relationship between
semigroups with context-free loop problem, and recent attempts to define
\textit{word hyperbolic} semigroups. A \textit{choice of representatives}
for a semigroup $S$ consists of a choice of generators $\sigma : X^+ \to S$ together
with a subset $R \subseteq X^+$ such that $R \sigma = S$. The choice of
representatives is called \textit{regular} if $X$ is a finite and $R$ is a
regular language. The \textit{multiplication table} for $S$ with respect
to $\sigma$ and $R$ is the language
$\lbrace u \# v \# w^R \mid (uv) \sigma = w \sigma \rbrace$
where $\#$ is a new symbol not in $X$, and $w^R$ denotes the word $w$
written backwards.

An interesting recent result of Gilman \cite{Gilman02} is that a finitely
generated group is word hyperbolic (in the sense of Gromov \cite{Gromov87})
if and only if it admits a regular choice of representatives with respect
to which the multiplication table is a context-free language. Duncan and
Gilman \cite{Duncan04} proposed that the latter condition might form a
suitable basis for a theory of \textit{word hyperbolic semigroups}.
Virtually free groups --- that is, groups with context-free word problem
--- form an elementary class of word hyperbolic groups. It transpires that
a more general relationship holds between semigroups with context-free loop
problem and word hyperbolic semigroups in the sense of Duncan and Gilman.
\begin{proposition}
Let $S$ be a semigroup with context-free loop problem. Let
$\sigma : X^+ \to S$ be any choice of generators for $S$. Then the
multiplication table of $S$ with respect to $\sigma$ and $X^+$ is
context-free. In particular, $S$ is word hyperbolic in the sense of
Duncan and Gilman \cite{Duncan04}.
\end{proposition}
\begin{proof}
Consider first the language
$$K_1 = \lbrace u \# v \# w \mid u, v, w \in \hat{X}^*, uvw \in L_\sigma(S) \rbrace$$
It is an easy exercise to verify that $K_1$ is a rational transduction of
$L_\sigma(S)$, and hence is context-free. Now by intersecting with a
regular set we see that the language
$$K_2 = K_1 \cap (X^* \# X^* \# \ol{X}^*)$$
is also context-free. Finally, applying the substitution $\ol{x} \to x$
for all $x \in X$ gives the multiplication table we require, so the latter
is also context-free.
\end{proof}

\section{Right Cancellative Monoids}\label{sec_canc}

In this section, we study the loop problem in the case of right
cancellative monoids. We begin with a proposition which tells us that,
in this case, the loop automaton is even more closely related to the loop
problem than in general.

\begin{proposition}\label{prop_minimality}
Let $\sigma : X^* \to M$ be a choice of generators for a right cancellative
monoid $M$. Then the loop automaton $\loopauto{\sigma}{M}$ is the
minimal automaton for the loop problem $L_\sigma(M)$.
\end{proposition}
\begin{proof}
By definition, $\loopauto{\sigma}{M}$ is an automaton accepting the loop
problem of $M$. Clearly there is a path from the identity vertex to every
vertex, and a path from each vertex to the identity vertex, so that the automaton
is trim.

For determinism, suppose $x \in X$ and $p, q, r \in M$ are such that there
are edges from $p$ to $q$ and from $p$ to $r$ both labelled $x$. Then we
have $q = p (x \sigma) = r$. On the other hand, suppose there are edges from $p$ to
$q$ and from $p$ to $r$ both labelled $\ol{x}$. Then $q (x \sigma) = p = r (x \sigma)$ so by right
cancellativity we deduce again that $q = r$. Thus, the automaton is
deterministic.

For minimality, it suffices to show that no two states in the loop
automaton are equivalent, that is, have the same cone.
To this end, let $p$ and $q$ be vertices in the loop automaton, that
is, elements of $M$, and choose a word $u \in X^*$ such that $u \sigma = p$.
Then there is a path from $1$ to $p$ labelled $u$, but no path from
$1$ to $q$ labelled $u$. But now there is an inverse path from $p$ to $1$
labelled $\ol{u}$, but no path from $q$ to $1$ labelled $\ol{u}$. Hence,
$\ol{u}$ lies in the cone of $p$ but not in the cone of $q$. We deduce that
all vertices have distinct cones, so that the automaton is minimal.
\end{proof}

Combining Proposition~\ref{prop_minimality} with Proposition~\ref{prop_syntrans}, we immediately obtain
the following description of the syntactic monoid of the loop problem of
a right cancellative monoid.
\begin{corollary}\label{cor_synloop}
Let $\sigma : X^* \to M$ be a choice of generators for a right cancellative
monoid $M$. Then the transition monoid of the loop automaton
$\loopauto{\sigma}{M}$ is the syntactic monoid of the loop problem
$L_\sigma(M)$.
\end{corollary}

We now describe an interaction between the loop problem and the classical
structural theory of right cancellative monoids. Let $M$ be a right
cancellative monoid. For each element $m \in M$, define the \textit{right translation map} $\rho_m : M \to M$ by $x \rho = xm$ for
all $x \in M$. Since the monoid is right cancellative, each map $\rho_m$ is
injective. Its \textit{(relational) inverse} is the map
$\rho_m^{-1} : M m \to M$ well-defined by $(xm) \rho = x$ for all $x \in M$.
Viewed as partial bijections on $M$, the maps of the form $\rho_m$ and
$\rho_m^{-1}$ generate an inverse monoid, which is called the \textit{inverse
hull} of $M$. Inverse hulls, were introduced by Rees \cite{Rees48}, who
used them to study embeddings of cancellative semigroups into groups. A
detailed study can be found in \cite[Sections~I.9 and I.10]{Clifford61}.

Now let $\tau : X^* \to M$ be a choice of generators for $M$.
It is readily verified from the definitions that for each $x \in X$, the
right translation map $\rho_{x \tau}$ is exactly the same as the transition
map $\sigma_x$ in the loop automaton $\loopauto{\tau}{M}$. Similarly, its relational inverse
$\rho_{x \tau}^{-1}$ is the transition map $\sigma_{\ol{x}}$ in the loop
automaton. Thus, we see that the transition monoid of the loop automaton
$\loopauto{\tau}{M}$ is exactly the inverse hull of the right cancellative
monoid $M$. Combining with Corollary~\ref{cor_synloop}, we obtain the
following relationship between the loop problem and the inverse hull of
a right cancellative monoid.

\begin{theorem}\label{thm_invhull}
Let $\tau : X^* \to M$ be a choice of generators for a right cancellative
monoid $M$. Then the syntactic monoid of the loop problem $L_\tau(M)$ is
the inverse hull of $M$. (Moreover, its action by partial maps on the
minimal automaton of $L_\tau(M)$ coincides with its action by partial
bijections on $M$.)
\end{theorem}
Note that the inverse hull of a group is easily seen to be isomorphic to
the group itself. Hence, in this case, we recover the ``folklore'' fact
that a group is the syntactic monoid of its own word problem.

In the case that $M$ is right cancellative, we can describe another
nice language-theoretic property of the loop problem. Recall that a
language $L$ is called \textit{deletion-closed} if whenever $x, y, w$ are
(possibly empty) words such that $w \in L$ and $xwy \in L$ we have
$xy \in L$. Note that a non-empty deletion-closed language will always
contain the empty word.
\begin{proposition}\label{prop_delclosed}
Any loop problem of a right cancellative monoid is deletion-closed.
\end{proposition}
\begin{proof}
Let $\sigma : X^* \to M$ be a choice of generators for a right cancellative
monoid $M$, and suppose $xwy, w \in L_\sigma(M)$. Then the loop automaton has vertices $p$
and $q$ such that there is a path $\pi_1$ from $1$ to $p$ labelled $x$, a
path $\pi_2$ from $p$ to $q$ labelled $w$ and a path $\pi_3$ from $q$ to $1$
labelled $y$. There is also a loop at $1$ labelled $w$, so by
Lemma~\ref{lemma_loopsdown} there is a loop $\pi_4$ at
$p$ labelled $w$. Now the paths $\pi_2$ and $\pi_4$ both start at $p$ and
have label $w$, and by Proposition~\ref{prop_minimality} the loop automaton
is deterministic, so $\pi_2$ and $\pi_4$ must end at the same vertex, that
is, $p = q$. But now $\pi_1 \pi_3$ is a loop at $1$ labelled $xy$, so
$xy \in L$ as required.
\end{proof}

Combining Propositions~\ref{prop_insertclosed} and~\ref{prop_delclosed}
we obtain the following corollary.
\begin{corollary}
Every inverse hull of a right cancellative monoid is the syntactic monoid
of an insertion-closed, deletion-closed language.
\end{corollary}

\begin{proposition}\label{prop_insdel}
Let $L \subseteq X^*$ be a formal language. Then the following are equivalent:
\begin{itemize}
\item[(i)] $L$ is insertion-closed and deletion-closed;
\item[(ii)] $L$ is the language of words representing the identity in its syntactic monoid;
\item[(iii)] $L$ is the language of words representing the identity in some monoid.
\end{itemize}
\end{proposition}
\begin{proof}
Suppose (i) holds. Observe that the map $w \mapsto [w]_L$ is a choice of generators for the
syntactic monoid $S(L)$. Now the fact that $L$ is insertion-closed and
deletion-closed means precisely that for any $w \in L$ and any $x,y \in X^*$,
we have $xwy \in L$ if and only if $xy \in L$. But this is true exactly if
$w$ is syntactically equivalent to the empty word, that is, if $[w]_L$ is
the identity element in $S(L)$. Hence, (ii) holds.

That (ii) implies (iii) is immediate, so suppose now that (iii) holds,
that is, that $L$ is the language of words representing the identity in
some monoid $M$ generated by $X$.
Clearly if $xy$ and $w$ both represent the identity in $M$ then $xwy$ also
represents the identity, so $L$ is insertion-closed. Conversely, if $xwy$
and $w$ represent the identity then $xy$ represents the identity, so $L$
is deletion-closed. Thus, (i) holds.
\end{proof}

Combining Theorem~\ref{thm_invhull} with Propositions~\ref{prop_delclosed}
and~\ref{prop_insdel} we obtain the following additional description of
the loop problem in a right cancellative monoid.

\begin{corollary}
Let $M$ be a right cancellative monoid. Then the loop problem of $M$ is
the language of words representing the identity in the inverse hull of $M$.
\end{corollary}

An application of the connections developed in this section to the study
of a class of inverse monoids will form part of a forthcoming paper of
J.~B.~Fountain and the author \cite{KambitesPolygraph}.

\section*{Acknowledgements}

This research was started while the author was at Carleton University
supported by the Leverhulme Trust, and completed at Universit\"at Kassel
with the support of a Marie Curie Intra-European Fellowship within the
6th European Community Framework Programme. The author would like to thank
John Fountain for some helpful discussions and Kirsty for all her support
and encouragement.

\bibliographystyle{plain}

\def\cprime{$'$} \def\cprime{$'$}

\end{document}